\documentclass{amsart}

\usepackage{amssymb}

\title[Perspective on symmetry groups]{A group-theorist's perspective on 
symmetry groups in
physics}
\author{Robert Arnott Wilson}
\date{First draft: 26th September 2020. This version: 20th December 2020}

\address{Queen Mary University of London}
\email{r.a.wilson@qmul.ac.uk}

\begin{document}
\begin{abstract}
There are many Lie groups used in 
physics, including the Lorentz group of special relativity,
the spin groups (relativistic and non-relativistic) and the gauge groups of quantum electrodynamics and the weak and strong nuclear forces.
Various grand unified theories use larger Lie groups in different attempts to unify some of these groups into something more fundamental.
There are also a number of finite symmetry groups that are related to the finite number of distinct elementary particle types.
I offer a group-theorist's perspective on these groups, and suggest some 
ways in which a deeper use of group theory
might in principle be useful. 
These suggestions include a number of options that seem not to be under active investigation at present. I leave open the question of whether 
they can be implemented in physical theories.
\end{abstract}
\maketitle

\section{Introduction}
\subsection{The status quo}
The most important Lie group in macroscopic physics is surely the Lorentz group $SO(3,1)$, that describes the coordinate transformations on spacetime that are necessary to incorporate the experimental fact that the speed of light in a vacuum is independent of the relative velocity of the source and the observer. This group is the cornerstone of the theory of special relativity \cite{SR}, and modern physics is almost inconceivable without it. The Lorentz group is also used as a foundation for the theory of general relativity \cite{GR1, GR2, GR3}, but it is also possible to use the larger group $SL(4,\mathbb R)$, that is the special linear group of degree $4$, consisting of all $4\times 4$ real matrices with determinant $1$. Either this group or the general linear group $GL(4,\mathbb R)$ is used in various attempts to quantise the theory of gravity \cite{GL4R1,GL4R2}. Of particular importance in all these cases is the representation of the group on the anti-symmetric square of spacetime. In special relativity, this representation holds the electromagnetic field, as a complex $3$-dimensional vector field. The Lorentz group acts in this representation as the complex orthogonal group $SO(3,\mathbb C)$. Again, it is possible to extend to $SL(4,\mathbb R)$ acting as the orthogonal group $SO(3,3)$.

 The Lie groups used in quantum physics \cite{Griffiths,Zee} are rather different. The most important is surely the spin group $SU(2)$, that describes the spin of an electron, and is necessary in order to incorporate the experimental fact that the electron has spin $1/2$. The spin group is isomorphic to a double cover of the orthogonal group $SO(3)$, and it can be extended \cite{Dirac} to the relativistic spin group $SL(2,\mathbb C)$ that is isomorphic to a double cover of the restricted Lorentz group $SO(3,1)^+$. Then there are the gauge groups $U(1)$ of electromagnetism, $SU(2)$ of the weak force 
 \cite{Weinberg} and $SU(3)$ of the strong force \cite{GellMann}. The weak gauge group is isomorphic to the spin group, but is certainly not equal to it. I therefore want to maintain a stricter than usual distinction between equality of groups, and isomorphism of groups. 

Various grand unified theories have attempted to combine these three gauge groups into a larger Lie group. Two important cases are the Georgi--Glashow model \cite{GG} built on $SU(5)$, and the Pati--Salam model \cite{PatiSalam} built on $SU(2) \times SU(2) \times SU(4)$. These theories are not generally accepted, since they appear to predict proton decay, which is essentially ruled out by experiment. However, they continue to be influential \cite{Baez}, and still appear to contain some important insights into the nature of elementary particles. Many larger Lie groups have also been used \cite{HH,E6,E8}, but all essentially suffer from the same problem, and/or other problems \cite{DG}.

It is also important to mention the finite symmetries, which in the standard model are not fundamental, but are derived from the Lie groups. For example, the weak gauge group has a $2$-dimensional complex representation, which is used to describe the discrete `symmetries'
 between certain pairs of particles, such as electron and neutrino, or proton and neutron, or up and down quark. This group also has a $3$-dimensional complex representation, which is used to describe the three intermediate vector bosons, that is the neutral $Z$ boson and the two charged $W^+$ and $W^-$ bosons. One goal of the grand unified theories is to unify all the fundamental bosons (force mediators) into a single representation, and similarly unify all the fundamental fermions. In the standard model, there are 13 fundamental bosons, and 45 fundamental fermions. The list of 45 fermions is almost certainly complete, but the list of $13$ bosons does not include any mediators for the force of gravity, so its completeness is perhaps still open to question.
 
 \subsection{Possible developments}
 One thing that immediately strikes a group-theorist faced with this list of groups is that macroscopic physics is mostly described by real and/or orthogonal Lie groups, while quantum physics is mostly described by complex and/or unitary Lie groups. A unified viewpoint must somehow reconcile these two types of groups.  The standard identification of the relativistic spin group with a double cover of $SO(3,1)^+$ goes some way towards this, but does not deal with the gauge groups. It is easy, of course, to identify $U(1)$ with $SO(2)$, and the weak $SU(2)$ is again isomorphic to a double cover of $SO(3)$, but the strong $SU(3)$ cannot be related to any real or orthogonal group in such a way. Conversely, the orthogonal group $SO(3,3)$ cannot be related in this way to a unitary group. 
 
 In other words, the standard model appears to present a broadly consistent approach across the range of theories from quantum electrodynamics and the weak force, to macroscopic electromagnetism and special relativity. But it cannot incorporate either the strong force or general relativity into the same group-theoretical system. There may of course be very good physical reasons for this, and a consistent unification of all four fundamental forces may not be possible. 
 On the other hand, \emph{if} a new approach to unification can be found, different from all those that have been tried before, \emph{then} it may be possible to make progress.
 
 In this paper I consider what mathematical changes might be necessary in order to produce such a unification. I pay relatively little attention to whether these mathematical changes are consistent with experimental reality, although I do point out a few experimental results that I think are relevant to deciding between different options. The two main questions to decide are (a) whether to try to transplant the macroscopic real/orthogonal groups to quantum physics, or the quantum complex/unitary groups to macroscopic physics, and (b) how to relate the finite groups to the Lie groups.
 
 \section{From large to small}
 \subsection{Assembling the gauge groups}
 I start by considering the possibility of using the macroscopic group $SL(4,\mathbb R)$ in quantum physics. This group is not itself an orthogonal group, but is isomorphic to the spin group $Spin(3,3)$, a double cover of the simple orthogonal group $SO(3,3)^+$. In other words, it can be used as a container for all the spin-type groups, including $SL(2,\mathbb C)$ and the electroweak gauge group $U(2)$. It does not, of course, contain $SU(3)$, so that some modification to quantum chromodynamics (QCD) would certainly be required in order to carry this programme through. Since direct (rather than indirect) evidence for QCD is hard to come by experimentally, this may not be too serious a problem. It is certainly possible to embed $SL(3,\mathbb R)$ in $SL(4,\mathbb R)$, and $SL(3,\mathbb R)$ is just a different real form of $SU(3)$, specifically the split real form rather than the compact real form. It is therefore plausible to model exactly the same calculations as are done in QCD, by the simple expedient of introducing a few judicious factors of $i$ into the calculations at strategic points.
 
 With this important proviso, therefore, $SL(4,\mathbb R)$ contains all the necessary subgroups. But are they in the right places? That is, is it possible to embed these subgroups in such a way that all the required relationships between them are correct? The first, obvious, point to make is that $SL(4,\mathbb R)$ does not contain the direct product of all the 
 specified subgroups. It is therefore necessary to 
 examine the assumption that these groups all commute with each other.  
 
 In fact, this assumption 
 is problematic even in the standard model. The gauge group $U(1)$ of quantum electrodynamics does not in the usual formalism 
 commute with the weak gauge group $SU(2)$, but only with a subgroup $U(1)$ of $SU(2)$. 
 In order to maintain 
 that the weak $SU(2)$ commutes with the strong $SU(3)$ it is necessary to implement a whole raft of mixing angles in the Cabibbo--Kobayashi--Maskawa (CKM) matrix \cite{Cabibbo,KM} and the Pontecorvo--Maki--Nakagawa--Sakata (PMNS) matrix \cite{Pontecorvo,MNS}. These mixing angles are therefore a strong hint that the groups \emph{in practice} (rather than in theory) do not commute with each other. If further proof is needed, the direct product of $SU(2)$ and $SU(3)$ (or any other real forms) cannot be embedded in a simple Lie group of dimension less than $24$, at which point there are too many gauge bosons for the model to be consistent with experiment.
 
 It is 
 surely necessary, at least, for every pair of the three groups to be disjoint. For the weak force, we need a copy of $SU(2)$, although in the standard model this is really the complex form $SL(2,\mathbb C)$, rather than a real form. For quantum electrodynamics (QED), we need a copy of $U(1)$, disjoint from $SU(2)$, and probably disjoint from $SL(2,\mathbb C)$. This is easy to arrange, since $U(1)$ defines a complex structure on the underlying $4$-dimensional real space, and there is a subgroup $SL(2,\mathbb C)$ of $SL(4,\mathbb R)$ preserving this complex structure. Both $U(1)$ and $SL(2,\mathbb C)$ are necessarily disjoint from $SL(3,\mathbb R)$, since every element of $SL(3,\mathbb R)$ fixes a non-zero vector, while only the identity elements of $U(1)$ and $SL(2,\mathbb C)$ fix any non-zero vector.
 
 \subsection{The Lorentz group}
 In macroscopic physics, the Lorentz group is a subgroup $SO(3,1)$ of $SL(4,\mathbb R)$. It is therefore not possible, in the scenario envisaged here, to identify the Lorentz group with the relativistic spin group $SL(2,\mathbb C)$. The latter group can still exist in the model, of course, and calculations can be done with it in exactly the same way as in the standard model. It is just that the relationship to macroscopic spacetime cannot work in quite the same way. 
 
 On the face of it, this looks like a proof that this proposed strategy cannot work. However, it is also possible that it shows why the relationship between quantum physics and classical physics is so problematical, and why the measurement problem is still unsolved. If the isomorphism between the Lorentz group and the central quotient of the relativistic spin group is a mathematical accident (since there is only one complex Lie algebra of rank $1$), rather than a fundamental physical principle, then the foundations of quantum mechanics appear in a
 different light \cite{remarks}.
 
 Indeed, the same problem 
 exists 
 already in non-relativistic quantum mechanics, in which the relevant groups are the spin group $SU(2)$ and the rotation group $SO(3)$. If these groups are related by being disjoint subgroups of $SL(4,\mathbb R)$, then they can generate the compact subgroup $SO(4)$ of $SL(4,\mathbb R)$. Since $SO(4)$ contains two normal subgroups isomorphic to $SU(2)$, there are two canonical quotient maps from $SO(4)$ onto $SO(3)$. Each of these canonical quotients restricts to a canonical isomorphism between the rotation group $SO(3)$ and the spin quotient $SO(3)$.
 Hence the requirement for the spin group to be canonically isomorphic to a double cover of the rotation group is still fulfilled. It is only that this canonical isomorphism cannot be mathematically represented as an equality, although it has exactly the same effect in physics as if it were an equality.
  
 Finally, it is worth remarking that the 
 usual assumption that the gauge groups commute with the Lorentz group is the conclusion of the Coleman--Mandula theorem \cite{ColemanMandula}. But there is no guarantee that the Coleman--Mandula theorem applies in the more general situation that I am considering here. The hypotheses of the theorem are quite restrictive, and only apply to a quite narrow class of potential theories. In particular, the hypotheses include the identification of the Lorentz group with a quotient of 
 the relativistic spin group $SL(2,\mathbb C)$, and this hypothesis does not hold in any potential model of the type being considered here.
 
 To conclude, it is possible to transplant the macroscopic groups into quantum mechanics, provided only that (a) QCD can be implemented with the split real form $SL(3,\mathbb R)$ rather than the compact real form $SU(3)$ of the gauge group, and (b) the foundations of quantum mechanics can be implemented with a canonical quotient (or pair of canonical quotients) from the spin group $SU(2)$ to the rotation group $SO(3)$, rather than assuming strict equality of groups. In group-theoretical terms, I believe this is a strong foundation on which to try to build a unified model.
 
 \subsection{The Clifford algebra} The relativistic spin group $Spin(3,1)$ can be constructed from a Clifford algebra \cite{Porteous} with either signature $(3,1)$ or $(1,3)$. In the former case, the algebra $Cl(3,1)$ is isomorphic as an algebra to the algebra of all real $4\times 4$ matrices. On the other hand, $Cl(1,3)$ requires $2\times 2$ quaternion matrices. In the standard model, these two algebras are combined into the algebra of all $4\times 4$ complex matrices. However well it may work in practice, this process is mathematically unnatural, and destroys the Clifford algebra property. 
  
 It is fairly clear that both $Cl(3,1)$ and $Cl(1,3)$ are required in the standard model. 
 Essentially, $Cl(3,1)$ is used for QED, and $Cl(1,3)$ for the weak interaction, although the `mixing' of the two forces 
 means that the modelling is rather more complicated than this. In the standard model each Clifford algebra is obtained from the other by multiplying the generators by $i$.
 But there are other, more natural, ways to combine the two algebras, for example by embedding them both in $Cl(3,3)$.  The structure of this Clifford algebra is discussed in detail in \cite{Furey,Cliffordsubs,INI14}.
  
 This Clifford algebra allows us to distinguish not only between $Spin(3,1)$ and $Spin(1,3)$, but also between these and yet another subgroup of $Spin(3,3)$ that is isomorphic to both of them. This is the subgroup 
 that maps to the subgroup $SO(3,\mathbb C)$ of $SO(3,3)$. 
 Distinguishing between these three isomorphic groups is likely to be useful in distinguishing between the left-handed and right-handed spinors, and potentially therefore explaining the chirality of the weak interaction. The standard model distinguishes only two different groups here. 
  
Since the concepts of left-handed and right-handed spins for massive fermions are not relativistically invariant, they are defined by the non-relativistic group $SO(3)\times SO(3)$ of $SO(3,3)$, corresponding to the subgroup $SO(4)$ of $Spin(3,3)$. This subgroup contains both a left-handed and a right-handed copy of $Spin(3)$, as well as a copy of $SO(3)$ that is the compact part of $SO(3,1)^+$.

 \subsection{Grand unified theories}
 If it is true that some of the unitary groups in the standard model of particle physics might be better replaced by orthogonal groups, then the same is likely to be true for some of the grand unified theories. In particular, the Georgi--Glashow model based on $SU(5)$ might be less in conflict with experiment if it were based on an orthogonal group $SO(5,\mathbb C)$ instead, or indeed on some real form thereof. The obvious real form to choose is $SO(2,3)$, since this exhibits a symmetry-breaking into $2+3$ that is a 
 feature of the standard model. 
 
 We have now reduced from a $24$-dimensional Lie group, with too many gauge bosons, to a $10$-dimensional one, with too few. The proposal made earlier to use the $15$-dimensional group $SL(4,\mathbb R)$ instead seems like the best compromise, as far as potential matching to experiment is concerned. An alternative is to consider replacing the unitary groups, not by orthogonal groups such as real or complex forms of  $SO(5)$ or $SO(6)$, but by finite groups of permutations on $5$ or $6$ letters.

 Indeed, if the theory is required to describe a finite number of distinct particles, then group theory is quite clear that a finite permutation group is required, not a Lie group. Moreover, finite groups are much more varied and subtle than Lie groups, and can describe many types of symmetries that Lie groups simply cannot describe. There are therefore many possibilities for generation symmetries, fermion-boson supersymmetries, and other types of symmetries that cannot be implemented in any grand unified theory based on Lie groups.

 \section{From small to large}
 \subsection{Two complex dimensions}
 Let us now consider instead the possibility that the standard model is essentially correct in all its details, including the universal use of unitary groups
 for the gauge groups, and look at the options for using unitary groups as symmetry groups on a large scale. The group $U(1)$ already is so used, essentially by multiplying the time coordinate of spacetime by $i$ so that the electromagnetic field has a 3-dimensional complex structure. This group does not act on spacetime, however, but only on the electromagnetic field.
 
 The spin group $SU(2)$ could in principle be implemented as a symmetry group of spacetime on any scale. It would impose a quaternionic structure on spacetime, such that time lies in the real part and $i,j,k$ represent three perpendicular space directions. As far as I can see, this is not inconsistent with any fundamental principles of physics. The natural metric on (non-relativistic) spacetime would then be Euclidean rather than Lorentzian (Minkowskian), but in a non-relativistic context that seems perfectly reasonable. 
 
 The extension to $SL(2,\mathbb C)$ is different, however, since it requires a choice of complex subalgebra in the quaternions, and therefore a choice of a direction in space. As long as this direction remains internal to the particle, and unobservable, there is no problem. But as soon as it becomes a macroscopic direction in space, it becomes observable and must be identified.
 In an inertial frame, there is no preferred direction in space, and this group cannot therefore have a sensible macroscopic meaning in an inertial frame. But we do not live in an inertial frame, and experiments are only rarely done in an inertial frame. In particular, expensive elementary particle experiments are never done in an inertial frame. 
 
 Therefore there are plenty of meaningful directions in macroscopic space that could potentially be used to define particular copies of $SL(2,\mathbb C)$ for use in particular contexts. In the case of the relativistic spin group, this macroscopic direction would appear to be local to the experiment, and unrelated to larger scales. However, experiments on the spin of entangled electrons can be on quite large scales, and it may be necessary to pay particular attention to how `local' is `local'.
 
 In the case of the weak gauge group, on the other hand, there is a distinct possibility that the observed symmetry-breaking of the weak force may indeed be related to a large-scale direction defined by the motion of the Earth's surface, to which the experiments are attached. If true, this would be a startling development, and call into question a number of basic assumptions about the nature of elementary particles. On the other hand, a direct experimental test of such a proposal  would be very challenging, and therefore I cannot see any good mathematical or physical reason for eliminating this possibility from consideration at this stage.
 
 \subsection{Three complex dimensions}
 The remaining group that we should consider is $SU(3)$. It must presumably act on a complex $3$-space. One obvious macroscopic complex $3$-space is the electromagnetic field. Is there any prospect that the symmetry group of the electromagnetic field might be better presented as $SU(3)$ rather than $SO(3,\mathbb C)$? This seems highly unlikely, given the phenomenal success of special relativity, but is it \emph{possible}? Indeed, a complex $3$-space in physics would normally be assumed automatically to be a unitary space. The fact that the electromagnetic field is treated instead as an orthogonal space should therefore make one sit up and ask questions. I merely ask the question. I do not profess to provide any answers.
 
 Another obvious complex $3$-space in electrodynamics is the momentum-current vector. Since the rest mass and charge are regarded as fixed and immutable, the symmetry group in classical electrodynamics is $SO(3)$. But is it possible that in general mass and charge can be converted into each other, and the real symmetry group is $SU(3)$? Experimental evidence is that rest mass is not conserved, but charge is, which makes it difficult to see how $SU(3)$ could act in this way. Nevertheless, is it possible 
 that there 
 is a more subtle way of interpreting the complex $3$-space, which might allow an $SU(3)$ symmetry group?
 Again, I merely ask the question, expecting the answer `no'.  
 
 Overall, I think my conclusion is that transplanting the complex/unitary groups into a macroscopic setting is less likely to be successful than transplanting the real/orthogonal groups in the opposite direction. The conflict with special relativity just seems insurmountable. But I cannot rule it out completely. The possibilities that general relativity (or some other theory of gravity)  is better described by $SU(4)$, or perhaps by $SU(3,1)$, than by $SL(4,\mathbb R)$, remain options, I believe.
  
  \section{The Klein correspondence}
  \subsection{Group isomorphisms and algebra isomorphisms}
  The isomorphism between the spin group $Spin(3)$ and the unitary group $SU(2)$ is part of a collection of isomorphisms between small
  spin groups and matrix groups, collectively known as the Klein correspondence.  This correspondence extends up to six dimensions, where the
  four different real forms of the spin group are:
  \begin{eqnarray}
  Spin(6) &\cong& SU(4)\cr
  Spin(5,1)&\cong& SL(2,\mathbb H)\cr
  Spin(4,2)&\cong& SU(2,2)\cr
  Spin(3,3)&\cong& SL(4,\mathbb R).
  \end{eqnarray}
  The five-dimensional versions are all symplectic groups:
  \begin{eqnarray}
  Spin(5) &\cong &Sp(2,\mathbb H)\cr
  Spin(4,1) &\cong& Sp(2,2,\mathbb R)\cr
  Spin(3,2) &\cong& Sp(4,\mathbb R)
  \end{eqnarray}
  The restrictions to 
  smaller dimensions are as follows:
 \begin{eqnarray}
 Spin(4) &\cong& SU(2)\times SU(2)\cr
 Spin(3,1)&\cong& SL(2,\mathbb C)\cr
 Spin(2,2)&\cong &SL(2,\mathbb R)\times SL(2,\mathbb R)\cr
 Spin(3)&\cong& SU(2)\cr
 Spin(2,1)&\cong& SL(2,\mathbb R)\cr
 Spin(2)&\cong & U(1)\cr
 Spin(1,1)&\cong& \mathbb R_{>0}^\times
 \end{eqnarray} 
  
  One can exhibit these isomorphisms in a number of different ways. From the point of view of group representation theory, one takes the group on
  the right hand side, acting on a $4$-dimensional space, and constructs the anti-symmetric square, so as to get a representation on a $6$-dimensional space. One then proves that this $6$-space has an invariant quadratic form of the given signature. 
  
  In the physics literature it is more usual to work
  instead from $6$ dimensions to $4$, which requires one to
  consider the whole $64$-dimensional Clifford algebra. 
  There are $7$ possible signatures on a real $6$-space, so $7$ distinct Clifford algebras, but among these there are only two distinct isomorphism types 
  of abstract algebra. The signatures $(6,0)$, $(5,1)$, $(2,4)$ and $(1,5)$ give $4\times 4$ quaternion matrix algebras, while the signatures $(4,2)$,
  $(3,3)$ and $(0,6)$ 
  lead to $8\times 8$ real matrix algebras. The different Clifford algebra structures in a given abstract algebra arise from
  different choices of generators. For example, if $A,B,C,D,E,F$ are six $8\times 8$ real matrices that square to $-1$ and anti-commute
  pairwise \cite{ManogueFairlie}, so that they generate $Cl(0,6)$, then $A,B,C,ABCD,ABCE,ABCF$ are generators for $Cl(3,3)$ and $A,B,DEF,CEF,CDF,CDE$ are
  generators for $Cl(4,2)$.
  
  In particular, the algebra of $8\times 8$ real matrices can be given Clifford algebra structures $Cl(4,2)$, $Cl(3,3)$ and $Cl(0,6)$, with
  respective spin groups $SU(2,2)$, $SL(4,\mathbb R)$ and $SU(4)$. These three groups are used in (a) the theory of twistors, (b) general
  relativity and (c) the Pati--Salam model. Of course, there is no reason to think that this mathematical coincidence can be translated into
  any physical relationship between these models. But \emph{if} there is any physical relationship between any two, or all three, of these theories \emph{then} this
  is a likely place in the mathematics to find it.
  
  \subsection{The Dirac algebra}
  Similarly, the Dirac algebra that plays a prominent role in quantum mechanics and 
  the standard model of particle physics is abstractly
  the algebra of $4\times 4$ complex matrices, which can be given Clifford algebra structures $Cl(4,1)$, $Cl(2,3)$ and $Cl(0,5)$. In all of these
  cases, the product of the five generators is the imaginary scalar $i$, so that 
 the usual interpretation is as a complex algebra with four generators,
  rather than a real algebra with five generators \cite{Cliffordsubs}. 
  
  The addition of a fifth generator corresponds mathematically to the embedding of $4$-dimensional spacetime into a larger space,
  similar to the construction of a de Sitter space or anti-de Sitter space from a $5$-space with signature $(4,1)$ or
  $(2,3)$ respectively. I leave aside all questions as to 
 any possible physical significance of
 this correspondence, 
 and merely remark on the fact of the 
  mathematical coincidence. Algebra is of no help in deciding such questions, which are all about the physical interpretation of the
  mathematics, not about the mathematics itself.
  
  In terms of the usual Dirac matrices $\gamma_0$, $\gamma_1$, $\gamma_2$ and $\gamma_3$ that are natural generators for $Cl(1,3)$,
  and $\gamma_5:=i\gamma_0\gamma_1\gamma_2\gamma_3$,
  suitable generators for other Clifford algebra structures can be chosen in various different ways, such as the following:
  \begin{eqnarray}
  Cl(4,1)&:& \gamma_0,i\gamma_1,i\gamma_2,i\gamma_3,i\gamma_5;\cr
  Cl(2,3)&:& \gamma_0,\gamma_5,\gamma_1,\gamma_2,\gamma_3;\cr
  Cl(0,5)&:& i\gamma_0,\gamma_1,\gamma_2,\gamma_3,i\gamma_5.
  \end{eqnarray}
  
  Extensions from $4\times 4$ complex matrices to $8\times 8$ real matrices can be obtained with an extra generator $j$ that squares to
  $-1$, commutes with
  $\gamma_0,\gamma_1,\gamma_2,\gamma_3$ and anti-commutes with $i$. On one hand, this algebra can be thought of as
  $Cl(1,3)$ with quaternion coefficients, while on the other hand it can be given various different structures as Clifford algebras of
  various different $6$-spaces which may or may not have good physical interpretations:
  \begin{eqnarray}
  Cl(0,6)&:&\gamma_1,\gamma_2,\gamma_3,i\gamma_0,j\gamma_0,k\gamma_0;\cr
  Cl(3,3)&:& \gamma_1,\gamma_2,\gamma_3,\gamma_0,j\gamma_5,k\gamma_5;\cr &&\gamma_1,\gamma_2,\gamma_3,\gamma_0,\gamma_5,k\gamma_5;\cr
  && \gamma_1,\gamma_2,\gamma_3,\gamma_5,j\gamma_5,k\gamma_5;\cr
&&i\gamma_1,i\gamma_2,i\gamma_3,i\gamma_0,j,k;\cr
Cl(4,2)&:& i\gamma_1,i\gamma_2,i\gamma_3,\gamma_0,j\gamma_0\gamma_5,k\gamma_0\gamma_5.
  \end{eqnarray}
  
  A useful mathematical criterion which may help in deciding between different choices of generators is that the product of
  all six generators (known as the pseudoscalar)
   is fixed by the relevant spin group. In the six cases listed above, this product is $i\gamma_5$, $\gamma_5$,
  $j\gamma_5$, $\gamma_0$, $\gamma_5$ and $i\gamma_5$ respectively. Another useful criterion could be based on what happens to
  the mass term in the Dirac equation under these various possibilities. In the usual Dirac equation, the mass term contains a factor of $i$,
  which commutes with the four Dirac matrices. Hence this can be generalised to $i$, $j$ and $k$, to give three mass terms for
  three generations of fermions. We may therefore want a pseudoscalar that is fixed under generation symmetry, and if so we may take
  either $i\gamma_5=\gamma_1\gamma_2\gamma_3\gamma_0$ or $\gamma_0$.
   
   A third criterion that may be useful is Lorentz invariance, which rules out $\gamma_0$. 
   Of the two remaining possibilities, only the first exhibits a clear generation symmetry. Moreover, the generation symmetry is
   associated only with the energy term, which seems a reasonable place to put something that determines three different
   masses, and therefore three different energies.
   
   \subsection{The compact case} 
   This leads us to consider the Clifford algebra $Cl(0,6)$ generated by $\gamma_1,\gamma_2,\gamma_3,i\gamma_0,j\gamma_0,k\gamma_0$
   as a likely candidate for generalising the Dirac algebra to three generations of fermions. The spin group is $SU(4)$. For a particular
   observer, spacetime clearly divides into space and time, so that we obtain a number of interesting subalgebras from this symmetry-breaking:
   \begin{eqnarray}
   Cl(0,3)&:& \gamma_1,\gamma_2,\gamma_3;\cr
   &&i\gamma_0,j\gamma_0,k\gamma_0;\cr
   Cl(1,3)&:& \gamma_1,\gamma_2,\gamma_3,\gamma_0;\cr
   &&\gamma_1\gamma_2\gamma_3,i\gamma_0,j\gamma_0,k\gamma_0;\cr
   Cl(3,0)&:&\gamma_0\gamma_1, \gamma_0\gamma_2, \gamma_0\gamma_3;\cr
   &&i\gamma_0\gamma_1\gamma_2\gamma_3,
j\gamma_0\gamma_1\gamma_2\gamma_3,
k\gamma_0\gamma_1\gamma_2\gamma_3.
   \end{eqnarray}
  The two copies of $Cl(1,3)$ are obtained by collapsing the first three or last three generators to their product. The first one defines the
  usual relativistic spin group, while the second defines another copy of $Spin(1,3)$ with generators
  \begin{eqnarray}
  i\gamma_0\gamma_1\gamma_2\gamma_3,
j\gamma_0\gamma_1\gamma_2\gamma_3,
k\gamma_0\gamma_1\gamma_2\gamma_3,&& i,j,k.
  \end{eqnarray}
  These two copies of $Spin(1,3)$ commute with each other, so that the second one can perhaps be taken to contain the gauge group $SU(2)$
  of the weak force, and the gauge group $U(1)$ of electrodynamics. It does not contain a copy of the gauge group $U(2)$
  of the standard model unification of these two forces, but this can be obtained by extending from $Spin(1,3)\cong SL(2,\mathbb C)$ to 
  $GL(2,\mathbb C)$, that is the subgroup of all invertible elements of the even part of the Clifford algebra. 
  Then $U(2)$ is the compact part of $GL(2,\mathbb C)$.
  
  It is not obvious, however, that this copy of $U(2)$, generated as a Lie group by $i,j,k$ for $SU(2)$ and $i\gamma_5$ for $U(1)$, behaves in the same way as the standard model gauge group. The latter, in practice, is implemented
  inside the Dirac algebra, that does not contain $j$ or $k$. Thus we have to restrict from $Cl(1,3)$ to $Cl(3,0)$ generated by
  \begin{eqnarray}
  \gamma_1\gamma_2\gamma_3,\gamma_0,\gamma_5.
  \end{eqnarray}
  This algebra is also a copy of the $2\times 2$ complex matrices, distinct from the even part of $Cl(1,3)$. Therefore it has its own subgroup $U(2)$,
  generated as a Lie group by $i\gamma_0$, $i\gamma_5$ and $\gamma_0\gamma_5$ for $SU(2)$, and $i$ for $U(1)$. 
  
  The algebra therefore contains two related copies of $U(2)$, one of which exhibits an unbroken symmetry between $i,j,k$, the other of which
  exhibits a broken symmetry between $i\gamma_0$, $i\gamma_5$ and $\gamma_0\gamma_5$. These are complexified with $i\gamma_5$ and $i$
  respectively, and in the standard model there is a  `mixing' between $i\gamma_5$ and $i$ that is the most conspicuous feature of
  the unification of QED and the weak interaction.
  
  Whatever the physical mechanism for this symmetry-breaking, the mathematical reason for it is the replacement of three generation labels $i,j,k$ 
  by a single mass label $i$. If we want to use the symmetric model for the three generations, then we cannot use $i$ as a mass label, but must use
  $i\gamma_5=\gamma_1\gamma_2\gamma_3\gamma_0$ instead. It is, of course, not obvious whether such a change can be incorporated into the
  standard model without breaking something important. Its effect on the Dirac equation is to extend the $SO(1,3)$ symmetry of energy-momentum,
  fixing the scalar mass term, to an $SO(1,4)$ symmetry of energy-momentum-mass \cite{Cliffordsubs}. The Clifford algebra $Cl(1,4)$ is generated by
  $\gamma_0,\gamma_1,\gamma_2,\gamma_3,\gamma_5$, and restricts to $Cl(3,1)$ generated by $i\gamma_0,i\gamma_1,i\gamma_2,i\gamma_3$, 
  or to $Cl(1,3)$ generated by $\gamma_0,\gamma_1,\gamma_2,\gamma_3$, 
 in both of which $i\gamma_5$ is the pseudoscalar. 
  
  It is worth also making some remarks about the spin copy of $Cl(1,3)$. The relativistic spin group embeds into $GL(2,\mathbb C)$ inside the
  even part of this algebra, in which the pseudoscalar element is again $i\gamma_5$. This contrasts with the scalar $i$ that is used in the standard model. 
  The distinction between left-handed and right-handed spin in the standard model is made by complex conjugation, that is, changing the sign of $i$.
  If instead we attempt to make this distinction by changing the sign of $i\gamma_5=\gamma_1\gamma_2\gamma_3\gamma_0$, then it
  becomes rather more obvious that this distinction really is a chirality.
  
  Indeed, the chirality has another mathematical realisation inside the subalgebra $Cl(0,3)$ generated by $\gamma_1,\gamma_2,\gamma_3$.
  This subalgebra is isomorphic to the direct sum of two quaternion algebras, obtained by projections with the orthogonal idempotents
  $1\pm\gamma_1\gamma_2\gamma_3$. The two algebras have generators
  \begin{eqnarray}
  \gamma_1+\gamma_2\gamma_3, \gamma_2+\gamma_3\gamma_1, \gamma_3+\gamma_1\gamma_2;\cr
  \gamma_1+\gamma_3\gamma_2, \gamma_2+\gamma_1\gamma_3, \gamma_3+\gamma_2\gamma_1;
  \end{eqnarray}
  from which the chirality is obvious.
  
  The other thing about the compact case that is particularly striking from a group-theoretical point of view
  is that the spin group is $SU(4)$, which suggests the possibility that it may
  contain a copy of $SU(3)$ that bears some relationship to the gauge group of the strong force. However, the Klein correspondence implies that
  any such $SU(3)$ imposes a unitary structure on the orthogonal space spanned by $\gamma_1,\gamma_2,\gamma_3, i\gamma_0,j\gamma_0,k\gamma_0$. The imaginary scalar would then link the three generations of fermions to three specific
  perpendicular directions in space. This is incompatible with the standard model, and incompatible with the Coleman--Mandula theorem.
  Hence this suggestion would appear to be ruled out.   
  
  \subsection{The split case}
  In the previous section I have shown how the compact case $Cl(0,6)$ gives rise to much of the structure of the standard model
  of particle physics, excluding the strong force, but including the three generations of fermions. 
  The pseudoscalar here is $i\gamma_5$, and is associated with a mass term that is independent of the
  three generations. Thus this entire structure is based on treating mass as an invariant. In order to incorporate the weak interaction,
  therefore, it was necessary to extend the groups beyond the spin groups,
  in order to allow the mass to change. 
  
  In the split case, on the other hand, the pseudoscalar is $\gamma_0$, associated with energy rather than mass. I ruled this case out
  on the basis that the model should be Lorentz-invariant. But if we regard energy as a more fundamental concept than mass, then
  Lorentz-invariance becomes a property of the observer's measurements (of mass), rather than a fundamental property of the underlying physics
  (of energy). So it is worth exploring the possibility of using the split case, based on the principle of conservation of energy,
  rather than the compact case, based on the principle of invariance of mass. In the end, of course, both Clifford algebras can be
  modelled with the same matrices, so there remains a possibility of translating between the two 
  viewpoints.
  
  The chosen generators for $Cl(3,3)$ can be written in one of the forms
  \begin{eqnarray}
  \gamma_1,\gamma_2,\gamma_3,&&i\gamma_0\gamma_1\gamma_2\gamma_3, j\gamma_0\gamma_1\gamma_2\gamma_3,
  k\gamma_0\gamma_1\gamma_2\gamma_3;\cr
  \gamma_1,\gamma_2,\gamma_3,&&\gamma_5,j\gamma_5,k\gamma_5.
  \end{eqnarray}
  The first exhibits the symmetry, while the second is easier to read.
  
  Breaking the symmetry between the three momentum coordinates and the three mass/generation coordinates gives us the following
  natural subalgebras:
  \begin{eqnarray}
  Cl(0,3) &:& \gamma_1,\gamma_2,\gamma_3;\cr
 && i\gamma_0,j\gamma_0,k\gamma_0;\cr
 &&\gamma_0\gamma_2\gamma_3, \gamma_0\gamma_3\gamma_1, \gamma_0\gamma_1\gamma_2;\cr
  Cl(0,4) &:& \gamma_1,\gamma_2,\gamma_3,\gamma_0\gamma_1\gamma_2\gamma_3;\cr
  Cl(3,0)&:& i\gamma_0\gamma_1\gamma_2\gamma_3, j\gamma_0\gamma_1\gamma_2\gamma_3,
  k\gamma_0\gamma_1\gamma_2\gamma_3;\cr
  Cl(4,0)&:&i\gamma_0\gamma_1\gamma_2\gamma_3, j\gamma_0\gamma_1\gamma_2\gamma_3,
  k\gamma_0\gamma_1\gamma_2\gamma_3;\gamma_1\gamma_2\gamma_3.
  \end{eqnarray}
  The first two copies of $Cl(0,3)$ and the
  copy of $Cl(3,0)$ have already been discussed in the context of the compact case.
  The third copy of $Cl(0,3)$ gives us a different way to split the spin group $SU(2)$ into two 
  pieces,
  given by the projections with $1\pm \gamma_0$. These are not chiral pieces, but may possibly have something to do with positive and negative charge,
  and/or particle/antiparticle pairs.
  
  The two algebras $Cl(4,0)$ and $Cl(0,4)$ give us copies of $Spin(4)\cong SU(2)\times SU(2)$
  that lie in the second and third copies of $Cl(0,3)$ respectively. The splitting of $Spin(4)$ 
  into two copies of $SU(2)$ is equivalent to the splitting of $Cl(0,3)$ into
  two copies of the quaternion algebra, and therefore contains the same physical information.
    
  The whole algebra $Cl(3,3)$ now defines a spin group isomorphic to $SL(4,\mathbb R)$. This is rather more difficult to interpret
  physically. The mathematical resemblance to the group $SL(4,\mathbb R)$ that describes the local coordinate transformations
  on spacetime between accelerating observers does not at first glance seem to correspond to a physical
  connection between the two theories.   But the appearance of $SL(4,\mathbb R)$ as a spin group does provide some mild 
  support for my earlier conclusion that this group may be useful
  in the theory of particle physics. 
  
  Conversely, the non-appearance of $SU(4)$ anywhere except as a spin group supports my other
  conclusion, that the unitary groups $SU(3)$ and $SU(4)$ are rather less likely to be useful in any realistic theory of
  macroscopic physics, including gravity. At any rate, the Klein correspondence does not support any such hypothetical applications.

 \section{From infinite to finite}
 \subsection{Finite symmetries of elementary particles}
  The second thing that strikes a group theorist as rather odd about the way that particle physicists use group theory, is the way that Lie groups are used to describe finite symmetries. If it works, of course, then there is no reason to change it. But it seems to me that it doesn't always work very well. There are certain places in the standard model where the number of experimentally observed particles of a given type does not match the dimension of the representation that is used to hold them. 
  
  The most obvious example concerns the kaons, which in theory live in two $2$-dimensional representations of $SU(2)$, but which experimentally are five (or possibly even six) distinct particles. It may well be that the mathematics that is used to deal with this situation describes kaons correctly. But this mathematics is not group theory, and the continued discovery of new 
 anomalies (i.e. disagreements between theory and experiment) for kaons \cite{CP,kaonanomaly,Kaon2} suggests that there may be a problem with the standard model in this area.

 \subsection{Finite permutation groups}
 The numbers of fundamental fermions of various types ($6$ leptons and $6\times 3$ quarks, consisting of $45$ individual Weyl spinors, $15$ in each of $3$ generations) suggest that the relevant group is closely related to the symmetric group $Sym(6)$ on $6$ letters, and the corresponding simple group $Alt(6)$. The latter is the most interesting and unusual of all the alternating groups \cite{FSG,Atlas}, since in addition to the extension to $Sym(6)$, and the double cover $2.Alt(6)$, that come from generic constructions applicable to all alternating groups, there is a further outer automorphism of $Sym(6)$ (a potential `supersymmetry'?) and a triple cover $3.Alt(6)$.
 
 The primitive permutation representations of $Alt(6)$ are on $6$, $15$, $10$, $15$ and $6$ points. The numbers $6$ and $10$ are the same as the numbers of anti-symmetric and symmetric rank $2$ tensors on $4$-dimensional spacetime, while $15$ is the dimension of the adjoint representation of $SL(4,\mathbb R)$, or of the compact real form $SU(4)$. More specifically, the Clifford algebra $Cl(3,3)$ or $Cl(6)$ splits into representations of dimensions $1+6+15+10+10+15+6+1$, such that, at least in the compact case, $Alt(6)$ can act as permutations on these representations. 
 
 Hence $Alt(6)$ can plausibly model finite symmetries on these various tensors and the Clifford algebra. The difference in signature between the compact real forms $SU(4)$ and $SO(6)$ on the one hand, and the split real forms $SL(4,\mathbb R)$ and $SO(3,3)$ on the other, suggests that there is an essential symmetry-breaking from $Alt(6)$ down to the subgroup that fixes a splitting of the $6$ letters into $3+3$. The further splitting into $1+2+3$ required by the standard model may arise from fixing a particular one of the $15$ fermions, which in the usual interpretation would be the right-handed electron. A different splitting as $2+1+3$ would permit an identification of the fixed fermion as 
 the left-handed neutrino instead.

There is a (unique) transitive permutation representation on 45 points, that might correspond to the standard model fundamental fermions.   Various subgroups of $Alt(6)$ can be used to describe the breaking of the symmetries between different particles. Of particular interest here are two conjugacy classes each of the groups $Alt(5)$, $Sym(4)$ and $Sym(3)$.  For example, 
 $Sym(3)$ acts as $1+2+3+3+6$ on one of the sets of $15$ points, which look 
 like the $15$ fermions in a single generation, 
 grouped 
 under $Alt(5)$ as $(2+3)+(1+3+6)$, 
 as in the Georgi--Glashow model. 

It is possible, of course, that the coincidences of the various small numbers mentioned here are just that, coincidences. But there is quite a limited choice of finite groups that can act on sets of these sizes, and \emph{if} there is a model of quantum physics based on finite groups, \emph{then} there is little choice but to use some combination of some of the groups I have suggested. 

\subsection{Finite matrix groups}
A further aspect of the finite perspective is the existence of finite representations, that is representations over finite fields \cite{ABC}. The group $Sym(6)$ is isomorphic to the $5$-dimensional orthogonal group $SO(5,2)$ over the field of two elements. This can be seen in the representation on all even-weight bit-strings of length $6$. There is a fixed bit-string, $111111$, so that the representation is not irreducible, but there is an irreducible representation on the even bit-strings modulo complementation. 

One of the subgroups $Alt(5)$ acts on this $4$-space as a $4$-dimensional orthogonal group, and can be used to describe the Georgi--Glashow allocation of $15$ fermions. For example, on the splitting A+BC+DEF of $6$ letters we can put the neutrino in AB, the left-handed electron in AC and the right-handed electron in BC, three colours of right-handed down quarks in AD, AE, AF and up quarks in DE, DF, EF, with the left-handed down quarks in CD, CD, CF and up quarks in BD, BE, BF. Then the splitting ABC+DEF reveals the leptons as a $2$-dimensional subspace of the $4$-space of even bit-strings modulo complementation, and the splitting BC+ADEF reveals a $3$-dimensional subspace containing all the right-handed particles. Thus the linear structure of the representation encodes certain fundamental physical properties, that are less easily visible from the Lie group perspective.

Moreover, one sees the transposition $(B,C)$ acting on the left-handed particles only, swapping electron with neutrino, and up with down quarks,
to create the weak doublets of the standard model. Similarly, the permutations of $D,E,F$ act to permute the three colours. It is noticeable that the right-handed up quarks appear to be different from the other quarks, in that they have two colours rather than one in this model. But colour confinement implies that this is equivalent to the standard model allocation of a single colour.

A different finite representation is provided by the isomorphism of the double cover $2.Alt(6)$ with the group $SL(2,9)$ of all $2\times 2$ matrices with determinant $1$ over the field of order $9$. This field is most easily obtained from the integers modulo $3$ (where $1+1=-1$) by adjoining a square root of $-1$, so forms a finite analogue of the complex numbers. Hence $SL(2,9)$ shares some of the properties of $SL(2,\mathbb C)$, that plays a fundamental role in the theory. Moreover, the subgroup $2.Alt(5)$ is isomorphic to $SU(2,5)$, and $2.Alt(4)$ is isomorphic to $SU(2,3)$, so there is a choice of various finite analogues of $SU(2)$ to use as required.

\subsection{Three generations}
The triple cover $3.Alt(6)$ has a $3$-dimensional representation over the field of order $4$, that is linear but not unitary. Here the field consists of elements $0,1,v,w$ such that $1+1=0$, $vw=1$ and $1+v+w=0$. The representation becomes unitary on restriction to the subgroup of index $10$ that fixes a given splitting of $6$ letters into $3+3$. This subgroup is isomorphic to a subgroup of index $2$ in $SU(3,2)$. 
The full group $SU(3,2)$ lies inside the larger group, sometimes called $3.M_{10}$, obtained by adjoining to $3.Alt(6)$ the outer automorphism that does not invert the scalars.  The group $3.Alt(6)$ acts on the $63$ non-zero vectors as $45+18$, so it may be possible to model the $45$ standard model fermions in this representation.

To see how this might work, take the following matrices as generators for $3.Alt(6)$: 
$$\begin{pmatrix} 0&1&0\cr 0&0&1\cr 1&0&0\end{pmatrix},
\begin{pmatrix}1&0&0\cr 0&w&0\cr 0&0&v\end{pmatrix},
\begin{pmatrix}0&1&0\cr 1&0&0\cr 0&0&1\end{pmatrix},
\begin{pmatrix}1&0&1\cr 0&1&1\cr 0&0&1\end{pmatrix}.
$$
These matrices map on to the permutations $(A,B,C)$, $(D,E,F)$,
$(B,C)(E,F)$ and $(A,D)(E,F)$ respectively.
Scalar multiplication implements the generation symmetry, so that ignoring the scalars gives back the usual fermions of a single generation, just as in the Georgi--Glashow model. With the same ordering of the letters as before, we can take the following types of particles:
\begin{eqnarray}
\nu=(1,0,0), & e_L=(0,1,0),& d_R=(1,1,0),\cr
e_R=(0,0,1), & u_R=(1,1,1), & u_L=(1,0,1),\quad  d_L=(0,1,1).
\end{eqnarray}
Colours and generations are obtained by multiplying the coordinates by $v$ and $w$ in suitable ways. The two rows give the $5+10$ splitting of the Georgi--Glashow model, acted on by a subgroup $Alt(5)\cong  
SL(2,4)$, acting on the first two coordinates. The standard model is obtained by restricting to a suitable copy of $SL(2,2)$.
 
 Two things are worth noting about this scheme. First, that the permutation $(D,E,F)$ that in the Georgi--Glashow model acts on colours only, here also acts on the generations of leptons. This is an unavoidable fact about this model, which cannot therefore implement a complete separation of the concepts of colour and generation as in the standard model. It is a matter of opinion whether one regards this as proof that this finite model cannot work, or instead as a hint as to how the three generations can be incorporated in a natural way into the standard model, rather than bolted on as an afterthought.
 
 Second, that the permutation $(B,C)$, that one might want to use to describe weak doublets, is odd, and therefore inverts the scalars $v$ and $w$ that are used for the triplet generation symmetry. To model a weak force that fixes the three generations of leptons it is necessary to use instead an even permutation, such as
 $(B,C)(E,F)$. Doing so creates an unavoidable mixing of the weak and strong forces, as they are described by the standard model. 
 The standard model does, of course, contain a mixing of the weak and strong forces, described by the CKM matrix. 
 
\section{From finite to infinite}
\subsection{Complex representations}
In order to relate these finite groups to Lie groups, we need to look at complex (and also real and quaternionic) representations. For this purpose it is useful to consider separately the representations of the double cover $2.Alt(6)$, the triple cover $3.Alt(6)$ and the simple group $Alt(6)$.
 There are representations of $3.Alt(6)$ inside $SU(3)$, that restrict to representations of $Alt(5)$ inside $SO(3)$, and may be related to the strong force with gauge group $SU(3)$. There are representations of $2.Alt(6)$ inside $Sp(4)$ inside $SU(4)$ that restrict to representations of $2.Alt(5)$ inside $SU(2)$ inside $SU(2)\times SU(2)$, and may be related to the weak force with gauge group $SU(2)$. Hence 
 several of the important ingredients of the standard model can be found in finite form in the representation theory of $Alt(6)$ and its covering groups.
 
\subsection{The double cover}
 Let us look first at the representations of $2.Alt(6)$. There are two $4$-dimensional complex representations, that are really $2$-dimensional quaternionic representations, and are therefore suitable for modelling spin. They both extend to $2.Sym(6)$, and are swapped by the outer automorphism of $Sym(6)$. Let us 
 consider the possibility of using one of them for the Dirac spinors.  Its anti-symmetric square breaks up as $1+5$, and its symmetric square is an irreducible $10$. To obtain the required structure of the complexified Clifford algebra as used in the standard model,  it is necessary to restrict to a subgroup $2.Alt(5)$, such that the representation breaks up as the sum of two $1$-dimensional quaternionic representations, so that the relevant Lie group is now $SU(2)\times SU(2)$. Then $5$ restricts as $1+4$, while $10$ restricts as $3a+3b+4$, so that the Clifford algebra structure appears as $1+4+3a+3b+4+1$ for the group $SO(4)$. 
 
 Of course, this is a complex representation, so that one can change the signature of the orthogonal group at will. Hence all the required structure of the Dirac spinors and the Clifford algebra arises from these representations. The additional feature that is present in this approach, that is not present in the standard model, is the finite symmetry group $2.Alt(5)$.
 
 The other $2$-dimensional quaternionic representation of $2.Alt(6)$ remains irreducible on restriction to $2.Alt(5)$. Its anti-symmetric square breaks up as $1+5$, and its symmetric square as $1+4+5$. These representations do not obviously appear in the standard model, but may possibly be of use in some extension to include quantum gravity. The $5$-dimensional representation is the spin $2$ representation of $SU(2)$, which certainly 
 chimes with the fact that in certain prospective models of quantum gravity, gravitons have spin $2$. A possible alternative interpretation, invoking the hypothetical `supersymmetry' between fermions and bosons, is that these representations describe massive bosons, in which case the $4$ perhaps represents the charged pions and kaons, and $1+5$ the neutral pions, kaons and eta meson, and possibly also the eta-prime meson.
 
 A third representation that may be of interest is the tensor product of the two $2$-dimensional quaternionic representations. This is a real $16$-dimensional representation that breaks up as $8a+8b$ for $Alt(6)$, and further as $3a+5+3b+5$ for $Alt(5)$. We shall see shortly how the $8$-dimensional representations can be regarded as copies of the adjoint representation of $SU(3)$. 
 
  \subsection{The triple cover}
 Now let us look at the representations of $3.Alt(6)$. There are four distinct $3$-dimensional unitary representations, coming in dual pairs $3A/3A^*$ and $3B/3B^*$. The representations $3A$ and $3B$ differ in the sign of $\sqrt{5}$. The tensor product of any of these representations with its dual (complex conjugate) breaks up as $1+8$, which permits the identification of these $8$-spaces with the adjoint representations of two distinct copies of $SU(3)$. Hence one or both of them can be used in the standard model to represent the strong force.
 
 The tensor product of $3A$ with $3B^*$ is an irreducible real $9$-dimensional representation, and is isomorphic to the tensor product of $3B$ with $3A^*$. This representation does not appear in the standard model, but if we restrict to $Alt(5)$, so that $3A$ and $3B^*$ become real orthogonal representations, then we obtain the spin $(1,1)$ representation of $SO(4)$. This is simply the compact version of the spin $(1,1)$ representation of $SO(3,1)$ that is used in general relativity for the Einstein tensor, among other things. So there is a potential application of this representation in a theory of quantum gravity.
 
 Of course, this is highly speculative, but the results of  
 this section show that the representation theory of the exceptional sextuple cover $6.Alt(6)$ of the alternating group $Alt(6)$ can give rise to 
 all the necessary mathematical structures needed for all the four fundamental forces of nature. I believe, therefore, that this group is a strong contender for a universal symmetry group for elementary particles. The double cover $2.Alt(6)$ contains all the structures necessary for the electroweak forces, and the triple cover $3.Alt(6)$ contains what is necessary for the strong force and the three generations of fermions. The group $Alt(6)$ itself contains all that is necessary for large-scale physics, in which individual elementary particle spins and colours are no longer relevant, but combine, in a way that must ultimately be described by the representation theory, to create a macroscopic concept of mass.

\section{Group algebras}

\subsection{General principles}
 The strategy outlined above for associating Lie groups with representations of finite groups is part of the standard theory of representations.
 To construct the group algebra, one takes a set
of linearly independent vectors, one for each element of the group. Then one just extends the group multiplication bilinearly, to get a multiplication
on the whole algebra. Thus one obtains an algebra whose dimension is equal to the order of the original group.

 Usually one works over the complex numbers, as this gives the simplest theory.
 In general, the complex group algebra of a finite group is isomorphic to a direct sum of complex matrix algebras, whereas the
 real group algebra will contain a mixture of real, complex and quaternion matrix algebras. 
 
 The general principle of representation theory is that the finite group acts on the algebra by multiplication: the usual convention in finite
 group representation theory is that it acts on the right. Then the matrix groups act on the left, to mix multiple copies of the same
 representation. Thus in the complex theory, the matrix groups commute with the finite group. The situation in the real theory is
 a little more subtle, since the real representations cannot distinguish between $i$ and $-i$, so that the complex structure is not
 well-defined. In this case, the complex matrix group does not commute with the finite group, unless the group is abelian.
 
 This paradigm is where we should expect to find much of the standard model, with finite symmetries of particles acting on the right,
 and continuous symmetries of gauge groups acting on the left. But one can also consider the finite group acting on both sides, or the
 matrix group acting on both sides. In the latter situation, the finite group has been lost, and one may hope to obtain instead
 a description of classical physics. In the former situation, the gauge groups have disappeared, and one may hope to obtain
 instead a \emph{finite} theory of quantum mechanics. 
 
 Then 
 the mixed case in the middle would have to be interpreted 
as measurements of elementary particles in a particular macroscopic environment.
In other words, the group algebra paradigm has the potential to address the measurement problem, and the emergence of
classical physics from quantum physics. Of course, it may not solve these problems, but at least it can address them.

  \subsection{An example}  
  For physical applications it is likely that the extra subtleties of the real theory will be useful. 
The group $2.Alt(4)$ discussed above is a good example that illustrates all the important points. 
 Since $2.Alt(4)$ has $24$ elements, we get a $24$-dimensional algebra this way. The complex irreducible
representations of the group have dimensions $1$, $1$, $1$ and $3$ (bosonic) and $2$, $2$ and $2$ (fermionic), so that the complex
group algebra has the structure
\begin{eqnarray}
\mathbb C + \mathbb C + \mathbb C + M_3(\mathbb C) + M_2(\mathbb C) + M_2(\mathbb C) + M_2(\mathbb C).
\end{eqnarray}
But over the real numbers the dimensions are $1,2,3$ (bosonic) and $4,4$ (fermionic), giving an algebra
\begin{eqnarray}
\mathbb R + \mathbb C + M_3(\mathbb R) + \mathbb H + M_2(\mathbb C).
\end{eqnarray}
It is now possible to convert this to a Lie group, by first taking out a real scale factor from each summand, and then removing the singular matrices,
so that we obtain the group
\begin{eqnarray}
U(1) \times SL(3,\mathbb R) \times SU(2) \times SL(2,\mathbb C).
\end{eqnarray}

This process has given us the three gauge groups of the standard model, except that $SU(3)$ is replaced by the split real form
$SL(3,\mathbb R)$, together with the relativistic spin group. 
It is possible, therefore, that the finite group $2.Alt(4)$ may contain
the essential key to the structure of the standard model. It may be related to the $4$-dimensional structure of spacetime in some way.
The double cover is clearly required in order to model fermions, and we note that both the spin group and the gauge group of the weak interaction
arise from fermionic representations. Odd permutations of $Sym(4)$ do not occur, and would appear to be unphysical.

\subsection{Symmetries and particles}
In order to obtain the gauge bosons, one needs to look at the symmetry group, that is the automorphism group of $2.Alt(4)$,
which is isomorphic to $Sym(4)$. This group contains $12$ even permutations, splitting into $1$ identity element, $3$ elements of
order $2$ and $8$ elements of order $3$. These perhaps correspond to the mediators in the standard model, that is $1$ photon,
$3$ intermediate vector bosons and $8$ gluons.
The odd permutations split into $6$ transpositions, and $6$ elements of order $4$. The obvious (though not necessarily correct)
physical interpretation is as fermions,
consisting of $6$ leptons and $6$ quarks. 

An alternative interpretation is that $Sym(4)$ is the symmetry group of a quantum of $4$-dimensional spacetime, which we might label by the
four letters $X,Y,Z,T$. Fixing one of the four letters, say $T$, splits the classes of elements as
\begin{eqnarray}
1+3+(2+6) + (3+3) + 6.
\end{eqnarray}
This might be interpreted as a splitting of the $8$ gluons into $2$ colourless and $6$ coloured, and a splitting of the $6$ leptons into
$3$ neutrinos and $3$ charged leptons. In other words, $T$ is related to charge, and $X,Y,Z$ to colours, in some way yet to be determined.
To relate the internal quantities $X,Y,Z,T$ to macroscopic concepts we must again look at the real representation theory, so that the
permutation representation splits as a scalar $X+Y+Z+T$ plus a $3$-dimensional representation spanned by $X-T$, $Y-T$ and $Z-T$.

These four vectors represent particles that we observe within the internal space spanned by $X$, $Y$, $Z$ and $T$. The particles
$X-T$, $Y-T$ and $Z-T$ correspond to the transpositions $(X,T)$, $(Y,T)$ and $(Z,T)$ which I suggested to interpret as leptons.
The particle $X+Y+Z+T$ might be related to the $4$-cycles, that is the quarks, but the $4$-cycles impose an ordering on $X,Y,Z$,
which has no significance for this particle. It is plausible that one cyclic ordering of $X,Y,Z$ corresponds to up quarks and the opposite
cyclic ordering to down quarks. These means that, once the quarks are chosen, one needs only three to eliminate all effects of the
ordering of $X,Y,Z,T$.
I suggest, therefore, that this particle is a baryon, consisting of three quarks. 
The simplest possibility is therefore that $X+Y+Z+T$ represents the proton, and $X-T$, $Y-T$, $Z-T$ represent the three generations
of electron, with $T$ representing the charge.
Or perhaps it is better to take $X+Y-Z-T$ and so on for the electrons, so that the four vectors are mutually orthogonal.

If we add up the particles obtained so far, with three protons in order
to cancel the charge of the three generations of electron, we obtain the vector $4X+4Y+4Z$. Perhaps this represents a collection of four
neutral baryons, such as neutrons? Or perhaps we need to take into account that when a neutron decays into an electron and a proton,
there is also a neutrino involved in the process? I have suggested that neutrinos correspond to the transpositions $(X,Y)$, $(X,Z)$
and $(Y,Z)$. One possibility is to relate $(X,Y)$ to $Z$, and so on, so that the neutrinos (or anti-neutrinos)
are the elementary particles $X$, $Y$ and $Z$. Or perhaps $X+Y-Z$ and so on.
So let us add these in, 
to get a total of $5X+5Y+5Z$, or five neutrons.

In other words, this finite model of elementary particles suggests that 
there may be some relationship between the sum of all six leptons and three protons on one hand, and five neutrons on the other.
This relationship holds in a 4-dimensional space over the real numbers, and therefore on a macroscopic scale. Taking out any dependence
on a direction in macroscopic space reduces the number of degrees of freedom from $4$ to $2$, and suggests that these are mass and charge.
I have already ensured that the charge is the same in both cases, so what about the mass?

Since measured neutrino masses are far too small to have any noticeable effect, the suggestion is essentially that
\begin{eqnarray}
\label{emutau}
m(e)+m(\mu)+m(\tau)+3m(p) &=& 5m(n).
\end{eqnarray}
Of these masses, the $\tau$ mass is the
least accurately known, by about three orders of magnitude. Hence
we can re-cast this equation as a prediction
of the $\tau$ mass.
I take
experimental values in $\mathrm{MeV}/c^2$ from CODATA 2014 \cite{CODATA} as
\begin{eqnarray}
m(e) &=& .5109989461(31)\cr
m(\mu) &=&105.6583745(24)\cr
m(p) &=& 938.2720813(58)\cr
m(n) &=& 939.5654133(58).
\end{eqnarray}

Then the
prediction is
\begin{eqnarray}
m(\tau) _p&=& 5\times939.5654133(58) -3\times 938.2720813(58) \cr & &- 105.6583745(24) -.5109989461(31)\cr
&=& 1776.84145(3).
\end{eqnarray}
This gives a prediction to $8$ significant figures, 
 well within current 
experimental 
uncertainty, given by 
\begin{eqnarray}
m(\tau)_e=1776.86(12). 
\end{eqnarray}

\subsection{More examples}
Perhaps what is most significant about the group $2.Alt(4)$ discussed above is that it is a finite subgroup of $SU(2)$ that acts irreducibly in the
$3$-dimensional (spin $1$) representation. There are exactly two other groups with these properties, namely $2.Sym(4)$ and $2.Alt(5)$.
Both are subgroups of $2.Alt(6)$ and have been discussed at some length above. It seems possible, therefore, that the group algebras of these two
groups may provide further insights.

Consider first the group algebra of $2.Sym(4)$. In fact there are two distinct groups of this shape, but only one of them is a subgroup of
$SU(2)$, and this is the one I shall consider here. The bosonic representations are all real, of dimensions $1,1,2,3$ and $3$, while the
fermionic representations are all quaternionic, of dimensions $1,1$ and $2$, so that the group algebra has the form
\begin{eqnarray}
2\mathbb R + M_2(\mathbb R) + 2M_3(\mathbb R) + 2\mathbb H + M_2(\mathbb H).
\end{eqnarray}
   The appearance of two copies of $3\times 3$ matrices here might permit one to be used for colours and the other for generations. This is
   perhaps an improvement on the version based on $2.Alt(4)$, which required a single copy of $3\times 3$ matrices to be used for both.
   The remaining parts of the algebra show an intriguing parallel between the real (or bosonic) part and the quaternionic (or fermionic) part.
   Perhaps, then, it is the quaternion structure that encodes the three generations, and not the $3\times 3$ matrix structure.

Converting to a Lie group by throwing away the real scale factors and the singular matrices gives us
\begin{eqnarray}
SL(2,\mathbb R) \times SL(3,\mathbb R) \times SL(3,\mathbb R) \times SU(2) \times SU(2) \times SL(2,\mathbb H).
\end{eqnarray}
Here we see $SL(2,\mathbb H)\cong Spin(5,1)$ as a different real form of $SU(4)\cong Spin(6)$, so that the fermionic part of
the algebra appears to be closely related to the Pati--Salam model. So there may be potential for building a model out of this algebra,
but it is far from obvious how to interpret the various pieces.

Turning now to $2.Alt(5)$, the bosonic representations are real of dimensions $1,3,3,4$ and $5$, and the fermionic representations
are quaternionic of dimensions $1,1,2$ and $3$. The Lie group in this case is
\begin{eqnarray}
&SL(3,\mathbb R) \times SL(3,\mathbb R) \times SL(4,\mathbb R) \times SL(5,\mathbb R) \times \cr
&SU(2) \times SU(2) \times SL(2,\mathbb H)
\times SL(3,\mathbb H).
\end{eqnarray}
Again, this group appears to be too big  to be useful, unless the factors are interpreted in quite a different way from the standard model.

\section{A worked example}
\subsection{The binary tetrahedral group}
In the previous section 
I have given some heuristic arguments to suggest that the group algebra of the binary tetrahedral group
may be closely related to the standard model of particle physics, and may therefore repay more detailed study.
First we need to understand the group itself.

A regular tetrahedron has four vertices, say $W,X,Y,Z$, such that every even permutation of the vertices can be achieved by
a rotation. For example, rotation about the axis through $W$ and the centre of $XYZ$ gives the $3$-cycles $(X,Y,Z)$ and $(X,Z,Y)$,
while rotation about the axis through the midpoints of $WX$ and $YZ$ gives $(W,X)(Y,Z)$. Thus the tetrahedral group is
isomorphic to the alternating group $Alt(4)$ on four letters. If we let 
\begin{eqnarray}
a&:=&(W,X)(Y,Z)\cr
b&:=&(X,Y,Z),
\end{eqnarray}
 then the group can be
defined abstractly by the relations 
\begin{eqnarray}
&a^2=b^3=(ab)^3=1.
\end{eqnarray}

The binary tetrahedral group is a double cover $2.Alt(4)$, obtained by adjoining signs in a particular way. As an abstract group, it is
defined by the relations
\begin{eqnarray}
a^4=b^3=(ab)^3=1, a^2b=ba^2.
\end{eqnarray}
The best way to see the group as a concrete (or real) group is to map $a$ and $b$ to the quaternions $i$ and $(-1+i+j+k)/2$
respectively. Then the representation $R$ as right-multiplications on the real $4$-space spanned by $1,i,j,k$ can be taken
as an explicit copy of the group. 

The elements and their orders are as follows:
\begin{eqnarray}
&\begin{array}{lll}
\mbox{Order} & \mbox{Abstract}&\mbox{Quaternion}\cr\hline
1 & e=a^4 & 1\cr
2 & a^2& -1\cr
3 & b, ab, ba, a^3ba, b^2, a^3b^2, b^2a^3, a^3b^2a& (-1\pm i\pm j\pm k)/2\cr
4 & a, b^2ab, bab^2, a^3, b^2a^3b, ba^3b^2& \pm i, \pm j, \pm k\cr
6 & a^2b, a^3b, ba^3, aba, a^2b^2, ab^2, b^2a, ab^2a & (1\pm i\pm j\pm k)/2\cr\hline
\end{array}
\end{eqnarray}

\subsection{The representations}
\label{representations}
There are five irreducible real representations, with characters (defined as the traces of the representing matrices) as follows:
\begin{eqnarray}
\begin{array}{l|rrrrr}
\mbox{Element} & e & a^2 & a & b & a^2b\cr
\mbox{Name}\backslash \mbox{Order} & 1 & 2 & 4 & 3 & 6\cr\hline
\mbox{R(eal)} & 4 & -4 & 0 & -2 & 2\cr
\mbox{S(pinor)} & 4 & -4 & 0 & 1 & -1\cr
\mbox{T(rivial)}&1&1&1&1&1\cr
\mbox{U(nitary)} & 2&2&2&-1&-1\cr
\mbox{V(ector)} & 3&3&-1&0&0\cr\hline
\end{array}
\end{eqnarray}
The representations $R$ and $S$ are fermionic, and $T$, $U$ and $V$ are bosonic. 
A useful representation is the transitive permutation representation $P$ on $8$ points,
obtained by adding signs to the four letters to get 
$\pm W$, $\pm X$, $\pm Y$ and $\pm Z$,
permuted as: 
\begin{eqnarray}
a&=& (W,X,-W,-X)(Y,-Z,-Y,Z),\cr
b&=& (X,Y,Z)(-X,-Y,-Z).
\end{eqnarray}
This  
representation  
breaks up into $S+T+V$. 
The representation $S$ is obtained on
a $4$-space spanned by vectors that can also be called $W,X,Y,Z$, such that the symbols $-W,-X,-Y,-Z$ are identified with the negatives of
the vectors $W,X,Y,Z$. This contrasts with $T+V$, in which $-W,-X,-Y,-Z$ are identified with $W,X,Y,Z$. 

The real group algebra splits (as a representation for the group) as
\begin{eqnarray}
T+U+3V+R+2S.
\end{eqnarray}
This representation is called the regular representation, and is equivalent to the transitive permutation representation on 
$24$ letters. 
The given decomposition can also be expressed as an algebra decomposition
\begin{eqnarray}
\mathbb R + \mathbb C + M_3(\mathbb R) + \mathbb H + M_2(\mathbb C).
\end{eqnarray}
Now we must be careful to note that the complex structure on $U$ is only defined up to complex conjugation: there are two ways
to express $U$ as a $1$-dimensional complex representation. The same applies to $S$, which has two different descriptions as
a $2$-dimensional complex representation. In the case of $R$, the automorphism group $SO(3)$ of the quaternions
gives us a $3$-parameter family of possible notations.

The space $3V$ is acted on by $2.Alt(4)$, by right-multiplications in the group algebra. This action commutes with the action of $M_3(\mathbb R)$
by left-multiplication, that mixes up the three copies of $V$ in all possible ways.  The action of $\mathbb C$ on $U$ is a complex scalar action,
that can be written either on the left or on the right, so also commutes with the action of the finite group.
The action of $\mathbb H$ by left-multiplication on itself likewise commutes with the right-action of $2.Alt(4)$.
In these 
three cases, therefore, we have a real scalar together with a `gauge group' which commutes with the action of $2.Alt(4)$, and is
respectively $SL_3(\mathbb R)$, 
$U(1)$ and $SU(2)$.

In the complex group algebra, $S$ can be given two distinct complex structures, such that the left-action by two copies of $M_2(\mathbb C)$
commutes with the right-action of $2.Alt(4)$. But in the real group algebra, the two copies are mixed together, and only the real part
$M_2(\mathbb R)$ commutes with $2.Alt(4)$.  In particular, and
counter-intuitively, the imaginary scalar matrices do not commute with the finite group.
This is, I believe, the reason why the standard model has to distinguish carefully between $i$ and $i\gamma_5$. 
The important point to consider now is how to relate the unbroken symmetry of $SU(2)$ acting on $R$, to the broken symmetry
of $SL_2(\mathbb R)$ acting on $2S$.

Before we do this, I remark that the group $SL(2,\mathbb C)$ that appears here as a summand of the group algebra has exactly the same
status in the algebra as all the other gauge groups. Logically, therefore, it should be the gauge group of one of the forces of nature.
There is only one force not yet accounted for, and that is gravity. Hence 
it may be worth bearing in mind
the possibility of interpreting $SL(2,\mathbb C)$ as the gauge group of quantum gravity.

\subsection{Classical physics}
To try to see if this makes sense, let us first consider how the matrix groups act by conjugation on the group algebra. Leaving aside the
three real and two complex scalars, what is left consists of the adjoint representations of $SL(1,\mathbb H)$, $SL(2,\mathbb C)$ and
$SL(3,\mathbb R)$. The adjoint representation of $SL(2,\mathbb C)$ is used in classical physics for the
electromagnetic field. 
This suggests using the adjoint representation of $SL(1,\mathbb H)$ for the Newtonian gravitational field.
This leaves the $8$-dimensional representation of $SL(3,\mathbb R)$ for general relativity.

Note, however, that the group algebra does not contain any Minkowski space representation of $SL(2,\mathbb C)$,
so does not contain anything that can be interpreted as special relativity. 
This is, of course, a problem for formulating general relativity, which is usually based on the Minkowski space representation
of $SO(3,1)$. But if we restrict both $SO(3,1)$ and $SL(3,\mathbb R)$ to $SO(3)$, then the adjoint representation is the sum of
the spin $1$ and spin $2$ representations, so by adding a scalar we can reconstruct something that has the same
structure as the Einstein tensor. Of course, there is no guarantee that a theory resembling general relativity can be
built on top of this piece of representation theory. 

An alternative might be to interpret $\mathbb H$ as an absolute spacetime, quantised as $R$. Then $M_2(\mathbb C)$ would represent
the modifications to absolute spacetime that arise from using 
electromagnetism to observe the universe. That is, the Minkowski space representation of $M_2(\mathbb C)$ represents
the spacetime that we observe via electromagnetic radiation, quantised as $2S$. If so, then we can perhaps regard $M_3(\mathbb R)$ as
representing the spacetime that we observe via gravitational radiation. Einstein's `curvature of spacetime' then arises from the fact that
gravitational spacetime does not match electromagnetic spacetime. Indeed, $M_3(\mathbb R)$ acts only on space, not on time, thus
leaving the model with a well-defined absolute timescale on which to measure gravitational time, or as we might call it, universal time.

Of course, there are many experiments that demonstrate that moving clocks run slow. Usually these are interpreted as saying that time
itself changes. But the group algebra model suggests a different interpretation, namely that it is only the relationship between
electromagnetic time and gravitational time that changes. In other words, it is the clock that changes, not the time. 
This is a remarkably common-sense interpretation, worthy, I believe, of serious consideration.

To interpret this picture in terms of classical physics, we have to restrict to $SO(3)$ acting on all three of these matrix algebras.
This splits the electromagnetic field into an electric field and a magnetic field. At the same time it splits the $8$ dimensions of
$SL(3,\mathbb R)$ into $3+5$, whereby the former represents the Newtonian gravitational field, and the latter represents 
Einstein's modification. 
There are also three scalars to worry about: one real scalar in $M_3(\mathbb R)$, representing gravitational mass, and one
complex scalar in $M_2(\mathbb C)$, representing electromagnetic mass and charge. Therefore, just as the model
distinguishes between gravitational time and electromagnetic time, it distinguishes between gravitational (or universal)
mass and electromagnetic
(or inertial) mass.

\subsection{The standard model of particle physics}
The standard model describes the interaction of the quantum world with the macroscopic world. To see this in the group
algebra, we must consider the algebra under left-multiplication by the gauge groups, rather than by conjugation. This
leaves the right-multiplication for the finite group. The dimension, $24$, is approximately equal to the number of unexplained real
parameters in the standard model, which suggests that they might be allocated in some way to the left-multiplications
by the gauge groups. Of course, it is unlikely to be as simple as that, as some of the parameters may be split between
different parts of the algebra.

The component $M_3(\mathbb R)$ looks likely to 
contain the $9$ masses of three generations of electron, up and down quarks. The mixing angles from the
CKM and PMNS matrices could plausibly go in $M_2(\mathbb C)$: they are normally allocated to $3\times 3$ complex matrices
rather than $2\times 2$ complex matrices, but that could be just a different representation of the same algebra. 
Then we have little choice but to allocate
the neutrino masses to $\mathbb H$, perhaps supplemented by the mass of the $Z$ boson (since this algebra represents the weak
interaction),
or possibly the Higgs boson. 
That leaves $T+U$ to
contain 
the strong coupling constant and the fine structure constant, and perhaps the mass of the $W$ bosons. 
Of course, this is only a suggested approximate allocation of parameters, and much more work will be required to
sort this out in detail.

The group algebra model, however, suggests that the matrix algebras should be interpreted as macroscopic properties, including
properties of the experimental apparatus, as well as the electromagnetic and gravitational background.
If so, then the number of genuine physical constants among these $24$ is only $1+1+3+1+2=8$.
Indeed, it is well-attested experimentally that a number of the $24$ parameters do indeed vary with the energy
scale. The group algebra suggests that some parameters may vary with the gravitational field. Since the 
experiments cannot be shielded from the gravitational field, this raises the alarming possibility
that they 
might depend on the
geographical location of the experiment, or its orientation in space. Or, worse, on the time of day, the time of year and the phase of the Moon.
Of course, if these variations are small, then they will be undetectable amongst the experimental noise.
And we still have eight genuinely universal constants to play with, if we can identify them precisely.

If the proposed allocation of parameters is approximately correct, then we must reduce the $9$ charged fermion masses to $3$,
presumably by removing the quark masses. Similarly the $4$ masses in $\mathbb H$ reduce to $1$, presumably by removing the
neutrino masses. The $8$ mixing angles reduce to $2$, but I have no suggestion how to do this. Finally, the two coupling constants
reduce to one, as one would expect in a unified theory. Of course, this is not precise enough to make any predictions, or to suggest
any experiments to test the group algebra model.

There is an additional possibility, that some of the parameters might depend on gravitational/acceleration properties like the
length of the day, or the length of the year, or, since we are probably only interested in dimensionless parameters, their ratio. If so,
then this ratio must have left its mark somewhere in the plethora of parameters. Searching for it systematically would be like
looking for a needle in a haystack, and would be very unlikely to produce any convincing evidence. That is, unless it had left its
mark in a very obvious place. 

\subsection{Looking in the obvious place}
Since I have allocated general relativity to the same algebra as the strong force, this `obvious place' must involve the
strong force.  The most obvious such place that I can think of is the proton/neutron mass ratio. 
Therefore my 
conjecture, based on the foregoing  
thorough analysis of the mathematical structure of the group algebra, is that the proton/neutron mass ratio may be
approximately related to the (average)
day/year time ratio $d$, via a currently non-existent quantum theory of gravity, that is closely related to quantum chromodynamics. 
Indeed, the model seems to allow the possibility that this hypothetical quantum theory may be independent of Newtonian gravity, and only implement a general relativistic correction.

A plausible formula of this type is
\begin{eqnarray}
m(n)/m(p) &\approx& 1 + d/2. 
\end{eqnarray}
Of course, this formula is not correct. There are all sorts of very good physical reasons why no such formula could ever be correct.
The best we can hope for is that it is a proxy for a much 
more complicated formula that might eventually emerge from
a fully-fledged quantum theory of gravity.
Substituting in some numerical values we have
\begin{eqnarray}
m(n)/m(p) \approx &939.565/938.272 &\approx 1.001378,\cr
1+d/2 \approx &1+ 1/730.52 &\approx 1.001369.
\end{eqnarray}
One coincidence like this does not justify a theory, of course. But it is a pretty amazing coincidence, if it is a coincidence.

\subsection{
The electron} 
Perhaps there is some more low-hanging fruit, for example the electron/proton mass ratio. For this, one needs
to look at the weak force in the broken symmetry form, that is the real subalgebra $M_2(\mathbb R)$ of $M_2(\mathbb C)$.
In electromagnetic terms, this must be interpreted as corresponding to the subgroup $SO(2,1)$ of the Lorentz group,
and therefore to a choice of a direction in space. In gravitational terms, on the other hand, it corresponds to a choice
of a subgroup $SL(2,\mathbb R)$ of $SL(3,\mathbb R)$, and therefore a choice of a \emph{different} direction in space.

There is therefore a single angle in macroscopic $3$-space that defines the relationship between these two directions.
If it is related to days and years, it must surely be the angle of tilt of the Earth's axis, say $\theta\approx 23.44^\circ$. 
This angle must manifest itself as a projection, hence in the form $\sin\theta$ or $\cos\theta$.  
It does not then take much imagination to 
propose the formula
\begin{eqnarray}
m(e)/m(p) &\approx& (d/2)\sin\theta.
\end{eqnarray}

I emphasise again that this formula is clearly not correct. It is at best a proxy for a much more complicated formula that might
exist at some time in the future, dependent on some hypothetical quantum theory of gravity.
Let us in any case substitute approximate values
\begin{eqnarray}
m(e)/m(p)\approx &.5109989/938.272&\approx .00054462,\cr
(d/2)\sin\theta \approx & .3978/730.52 &\approx .00054453.
\end{eqnarray}
  Again we see a very striking coincidence, if it is a coincidence.
  
  \subsection{Quantum mechanics}
  I have looked at the group algebra in classical form, with the matrix groups acting on both sides, and seen that the algebra 
  contains just the right ingredients to mimic the classical theories of electrodynamics and gravity, in the form of general relativity.
  I have looked at the group algebra in mixed form, with the matrix group acting on one side, and the finite group acting on the
  other, and seen that the algebra contains just the right ingredients to mimic the standard model of particle physics.
  I have looked at the relationship between the two, and have provided evidence that the model accurately describes 
  at least one or two aspects of the
  emergence of classical physics from quantum physics.
  
  Now let us look at the algebra in finite form, that is with the finite group acting on both sides, and without the matrix groups.
  In this form the algebra contains no mass, no energy, no momentum. Indeed, it contains no continuous variables at all, only
  discrete quantum numbers. The relationship between this and the standard model, therefore, is exactly the same as the
  relationship between the standard model and classical physics. The wave-function and all the continuous variables that are used in
  the theory of quantum mechanics lie in the mixed form of the group algebra. In particular, in the group algebra model 
  the wavefunction is not an intrinsic property of an elementary particle, but a description of the average behaviour of an
  elementary particle placed inside a particular macroscopic environment. 
  
  Collapse of the wave-function, therefore, can be viewed not as a physical process,
  but as a difference in point of view. That is, do we consider the elementary particles interacting in isolation, or in
  a wider context?
  In other words, the group algebra model 
  suggests
  a new way to look at the measurement problem. 
  The probabilities and other continuous variables that appear in the standard theory all arise from putting the elementary
  particles into a macroscopic context. 
  
  Decoherence then occurs when the macroscopic context changes, or becomes impossible
  to define precisely. 
  In the usual experiments with electron spin or photon polarisation, 
  the standard approach is to assume that 
only the electromagnetic context
  is relevant, 
  and not the gravitational context.  But since the experiment cannot be shielded from any change in the
  gravitational environment, it would seem to be 
   both theoretically and experimentally impossible to rule out the possibility that quantum mechanical experiments of this kind
  are affected by the gravitational environment.  
  If so, then our knowledge of the gravitational environment on a large scale
  enables us to define the gravitational environment at both ends of a macroscopic experiment \emph{without} the particles
  involved communicating with each other directly.
    This may be enough to explain the illusion of superluminal transfer of information in experiments with
  entangled particles.
  
    One class of experiments in which the gravitational context does appear to play a role is in the
  detection of neutrino oscillations. Since neutrinos only interact with the electromagnetic environment via the weak interaction,
  a neutrino in flight is not affected by changes in the electromagnetic environment. 
  But in large-scale experiments the gravitational context clearly changes from one end of the
  experiment to the other. Similar considerations may apply to other puzzling experiments, such as those that are usually interpreted as
  CP-violation for neutral kaons. The gravitational context may not change much over a length of $57$ feet, but it is enough to be
  detected in the kaon decay experiments.
  
 \section{Conclusions}
 In this paper I have looked at three different ways in which the group theory that is used in fundamental physics might in principle be developed to a more sophisticated level. I have not tried to address in detail the question as to whether these suggested developments can in fact be applied to physics, although I have tried to evaluate at a basic philosophical level the relative merits of the different possible approaches. The three approaches can be summarised as
 \begin{enumerate}
 \item try to use the Lie groups from general relativity in particle physics;
 \item try to use the Lie groups from particle physics in (quantum) gravity;
 \item use finite groups in particle physics, and derive the Lie groups from representations of the finite groups.
 \end{enumerate}
 My initial conclusions are that the second option is unlikely to be viable, but that both the first and the last options show considerable promise.
 The first option looks likely to offer minor changes to the standard model that would permit a more thorough unification of the electromagnetic, weak and strong forces. 
 The last option offers a tantalising glimpse of a possible quantum theory of gravity, that is potentially quite different from existing approaches
\cite{string,SUSY,Rovelli,verlinde,hossenfelder}.
 
 Finally, I have considered a particular example of this last option that I believe offers the best hope for progress, and looked at
 what it might have to say about the big problems of emergence of classical physics from quantum physics, and the measurement problem.
 It is early days, but
 I believe 
the omens are good.

\end{document}